\documentclass[12pt]{amsart}

\usepackage{accents}
\usepackage{appendix}
\usepackage{amsfonts}
\usepackage{amsmath}
\usepackage{amssymb}	
\usepackage{amsthm,bm}
\usepackage{array,booktabs,multirow}
\usepackage{braket}
\usepackage{centernot}
\usepackage{cite}
\usepackage{comment}
\usepackage{dsfont}
\usepackage{mathalfa}
\usepackage[shortlabels]{enumitem}
\usepackage{etoolbox}
\usepackage{float}
\usepackage[hang, flushmargin]{footmisc}
\usepackage{latexsym}
\usepackage{lipsum}
\usepackage{needspace}
\usepackage{tikz}
\usepackage{hyperref}
\usetikzlibrary{matrix,arrows}
\usepackage{diagbox}
\usepackage{adjustbox}

\makeatletter
\def\namedlabel#1#2{\begingroup
   \def\@currentlabel{#2}%
   \label{#1}\endgroup
}
\makeatother

\theoremstyle{plain}
\newtheorem{thm}{Theorem}[section]
\newtheorem{cor}[thm]{Corollary}
\newtheorem{lem}[thm]{Lemma}
\newtheorem{prop}[thm]{Proposition}

\theoremstyle{definition}

\newtheorem{defn}[thm]{Definition}
\newtheorem{exam}[thm]{Example}
\newtheorem{prob}[thm]{Problems}

\theoremstyle{remark}

\setlist[enumerate,1]{leftmargin=2em}

\def\C{\mathbb C}

\def\H{\mathcal H}

\def\Pd{\mathcal P_d(\mathbb R)}
\def\R{\mathbb R}
\def\T{\mathbf T}

\def\U{U(\mathfrak{sl}_2)}

\newcommand{\floor}[1]{\left\lfloor #1 \right\rfloor}

\title[Treating dual Hahn polynomials as Racah polynomials]{A way to treat dual Hahn polynomials as Racah polynomials via the theory of Leonard pairs
\\[2ex] 
{\normalfont \normalsize \ttfamily  Dedicated to Paul Terwilliger on his 70th birthday}    
}

\author{Hau-Wen Huang}

\address[H.-W. Huang]{
Department of Mathematics\\
National Central University\\
Chung-Li 32001 Taiwan}
\email{hauwenh@math.ncu.edu.tw}

\begin{document}

\begin{abstract}
The dual Hahn polynomials $\{u_i(x)\}_{i=0}^d$ are a family of discrete orthogonal polynomials involving two real parameters $r$ and $s$. 
Let $L,L^*$ denote the corresponding Leonard pair.
Assume that $r\not=0$ and $r+s=0$. We show that $L,(L^*+\frac{r-d}{2})^{2}$ is a Leonard pair.
According to the theory of Leonard pairs, the polynomials $\{u_i(x)\}_{i=0}^d$ are not only the dual Hahn polynomials but also the Racah polynomials with respect to the same inner product.
\end{abstract}

\maketitle


{\footnotesize{\bf Keywords:} dual Hahn polynomials, Racah polynomials, Leonard pairs}

{\footnotesize{\bf MSC2020:} 05E30, 16S30, 33C45}

\allowdisplaybreaks

\section{Introduction}\label{sec:intro}
Suppose that $D\geq 3$ is an integer. The $D$-dimensional hypercube $H(D,2)$ is a bipartite distance-regular graph whose vertex set is $\{0,1\}^D$ and two vertices $x,y\in  \{0,1\}^D$ are adjacent if and only if $x$ and $y$ differ in exactly one coordinate. 
A notable result of the paper \cite{hypercube2002} is to give 
an algebra homomorphism from the universal enveloping algebra $\U$ of $\mathfrak{sl}_2$ onto the Terwilliger algebra of $H(D,2)$  \cite{TerAlgebraI,TerAlgebraII,TerAlgebraIII}. 
Since about that time, there have been a significant number  of papers on the Terwilliger algebras and Lie algebras, for instance \cite{halved:2023,halved:2024,Huang:CG&Hamming,
Huang:CG&Grassmann,Huang:CG&Johnson,Johnson:2007,dualpolar:2022,Doob2023,Doob2021,Johnson:2021,Hamming:2015,bipart&so6,Wen2024}.

Let $\frac{1}{2}H(D,2)$ denote the halved graph of $H(D,2)$. 
This paper is motivated by recent studies on the Terwilliger algebra of $\frac{1}{2}H(D,2)$.
The work \cite{halved:2023} pointed out that an algebra homomorphism $\natural$ from the universal (dual) Hahn algebra $\mathcal H$ \cite{Hahn:2019-1,Hahn:2019-2,Johnson:2021,Huang:CG&Johnson} into $\U$ is hidden behind the Terwilliger algebra of $\frac{1}{2}H(D,2)$. 
As an association scheme, there are additional $P$- and $Q$-polynomial structures for $\frac{1}{2}H(D,2)$ with odd $D$ \cite{BannaiIto1984,DRG_book1989}. 
A subsequent paper \cite{halved:2024} implies that another $Q$-polynomial structure is related to an algebra homomorphism $\sharp$ from the universal Racah algebra $\Re$ \cite{SH:2017-1,SH:2019-1,SH:2020,R&LD2014,gvz2013,gvz2014} into $\U$. 
By \cite[Theorem 1.5]{halved:2023} and \cite[Theorem 7.8]{halved:2024} the image of $\sharp$ is properly contained in that of $\natural$. Consequently, the algebra homomorphism $\natural$ cannot be factored into a composition of an algebra homomorphism $\H\to \Re$ followed by $\sharp$. 
A reverse question is whether there exists an algebra homomorphism $\Re\to \H$ such that the following diagram is commutative \cite[Problem 6.1]{Wen2024}:

\begin{table}[H]
\centering
\begin{tikzpicture}
\matrix(m)[matrix of math nodes,
row sep=3em, column sep=2em,
text height=1.5ex, text depth=0.25ex]
{
\Re
&
&\H
\\
&\U\\
};
\path[->,font=\scriptsize,>=angle 90]
(m-1-1) edge (m-1-3)
(m-1-3) edge[bend left] node[right] {$\natural$} (m-2-2)
(m-1-1) edge[bend right] node[left] {$\sharp$} (m-2-2);
\end{tikzpicture}
\end{table}

\noindent This algebraic problem remains unresolved.
Nevertheless, some results are achieved from the perspective of discrete orthogonal polynomials.

To state the results, we begin by recalling the definition of discrete orthogonal polynomials. 
Fix an integer $d\geq 0$. 
Let $x$ denote an indeterminate. 
Let $\Pd$ denote the vector space over the real number field $\R$ consisting of all $f(x)\in \R[x]$ with $\deg f\leq d$.
Pick $d+1$ distinct real scalars $\{\theta_i\}_{i=0}^d$ and $d+1$ positive real scalars $\{k_i^*\}_{i=0}^d$. 
There is an inner product $\langle \, ,\,\rangle$ on $\Pd$ defined by 
$$
\langle f, g\rangle=
\sum_{i=0}^d 
f(\theta_i) g(\theta_i) k_i^*
\qquad 
\hbox{for all $f,g\in \Pd$}.
$$
A sequence $\{f_i(x)\}_{i=0}^d$ in $\Pd$ is called a family of {\it discrete orthogonal polynomials} on the points $\{\theta_i\}_{i=0}^d$ with weights $\{k_i^*\}_{i=0}^d$ provided that $\deg f_i=i$ for all $i=0,1,\ldots,d$ and 
$$
\langle f_i ,f_j \rangle
=
\left\{
\begin{array}{ll}
0 \qquad 
&\hbox{if $i\not=j$},
\\
>0
\qquad 
&\hbox{if $i=j$}
\end{array}
\right.
\qquad 
(0\leq i,j\leq d).
$$

Set $r$ and $s$ as two real parameters in the interval $(-1,\infty)$. Recall the Pochhammer symbol 
$$
(x)_i=x(x+1)\cdots (x+i-1)
\qquad 
\hbox{for any integer $i\geq 0$}.
$$
Note that a vacuous product is interpreted as $1$. 
We will use the following parameters associated with $r$ and $s$:
\begin{align}
\theta_i
&=(d-i)(d-i+r+s+1)
\qquad 
(0\leq i\leq d),
\label{thetai}
\\
\theta_i^*
&=i
\qquad 
(0\leq i\leq d),
\label{thetai*}
\\
b_i&=(d-i)(d-i+s)
\qquad 
(0\leq i\leq d-1),
\label{bi}
\\
c_i&=i(i+r)
\qquad 
(1\leq i\leq d),
\label{ci}
\\
a_i&=
\left\{
\begin{array}{ll}
\theta_0-b_0
\qquad &\hbox{if $i=0$},
\\
\theta_0-b_i-c_i
\qquad &\hbox{if $i=1,2,\ldots,d-1$},
\\
\theta_0-c_d
\qquad &\hbox{if $i=d$}
\end{array}
\right.
\qquad 
(0\leq i\leq d),
\label{ai}
\\
k_i&=\frac{b_0b_1\cdots b_{i-1}}{c_1c_2\cdots c_i}
\qquad 
(0\leq i\leq d),
\label{ki-1}
\\
\nu&=
\frac{(\theta_0-\theta_1)(\theta_0-\theta_2)\cdots (\theta_0-\theta_d)}{c_1 c_2\cdots c_d},
\label{nu-1}
\\
b_i^*&=
\frac{(d-i)(i-d-s)(2d-2i+r+s+2)_i}
{(2d-2i+r+s)_{i+1}}
\qquad (0\leq i\leq d-1),
\label{bi*}
\\
c_i^*&=
\frac{i(i-d-r-1)(d-i+r+s+1)_{d-i}}
{(d-i+r+s+2)_{d-i+1}}
\qquad 
(1\leq i\leq d),
\label{ci*}
\\
a_i^*&=
\left\{
\begin{array}{ll}
\theta^*_0-b^*_0
\qquad &\hbox{if $i=0$},
\\
\theta^*_0-b^*_i-c^*_i
\qquad &\hbox{if $i=1,2,\ldots,d-1$},
\\
\theta^*_0-c^*_d
\qquad &\hbox{if $i=d$}
\end{array}
\right.
\qquad 
(0\leq i\leq d),
\label{ai*}
\\
\label{ki*-1}
k_i^*&=
\frac{b^*_0 b^*_1\cdots b^*_{i-1}}{c^*_1 c^*_2\cdots c^*_i}
\qquad (0\leq i\leq d).
\end{align}
By (\ref{thetai}) and since $r,s\in (-1,\infty)$ the scalars $\{\theta_i\}_{i=0}^d$ are mutually distinct. 
Substituting (\ref{bi}) and (\ref{ci}) into (\ref{ki-1}) yields that 
\begin{gather}
\label{ki-2}
k_i={d\choose i} \frac{(d-i+s+1)_i}{(r+1)_i}
\qquad 
(0\leq i\leq d).
\end{gather}
Substituting (\ref{bi*}) and (\ref{ci*}) into (\ref{ki*-1}) yields that 
\begin{gather}
\label{ki*-2}
k_i^*={d\choose i}
\frac{(-d-s)_i(d+r+s+1)_d}
{(-d-r)_i(2d-2i+r+s+2)_i(d-i+r+s+1)_{d-i}}
\qquad 
(0\leq i\leq d).
\end{gather}
Substituting (\ref{thetai}) and (\ref{ci}) into (\ref{nu-1}) yields that 
\begin{gather}
\label{nu-2}
\nu=\frac{(d+r+s+1)_d}{(r+1)_d}.
\end{gather}
Since $r,s\in(-1,\infty)$ the scalars $\{k_i\}_{i=0}^d$, $\{k_i^*\}_{i=0}^d$ and $\nu$ are positive.

Recall that the hypergeometric function is denoted and defined by
$$
{}_mF_n
\left(\genfrac..{0pt}{}{\alpha_1,\alpha_2,\ldots,\alpha_m}{\beta_1,\beta_2,\ldots,\beta_n}\,\bigg|\,x\right)
=\sum_{i=0}^\infty 
\frac{(\alpha_1)_i(\alpha_2)_i\cdots (\alpha_m)_i}
{(\beta_1)_i(\beta_2)_i\cdots (\beta_n)_i} \frac{x^i}{i!}
$$
where $m$ and $n$ are positive integers and $\{\alpha_i\}_{i=1}^m$ and $\{\beta_i\}_{i=1}^n$ are scalars. If a numerator parameter $\alpha_i$ is a non-positive integer, then the series is finite.

\begin{lem}
\label{lem:3F2}
There are unique polynomials $\{u_i(x)\}_{i=0}^d$ with $ u_i\in \Pd$ for all $i=0,1,\ldots,d$ satisfying 
\begin{align*}
u_i(\theta_j)
&={}_3F_2
\left(\genfrac..{0pt}{}{-i,-j, j-r-s-2d-1}{-s-d,-d}\,\bigg|\,1\right)
\qquad 
(0\leq i,j\leq d).
\end{align*}
\end{lem}
\begin{proof}
Since the scalars $\{\theta_i\}_{i=0}^d$ are distinct, every polynomial $f\in \Pd$ is uniquely determined by the values of $f$ at $\{\theta_i\}_{i=0}^d$. The lemma follows.
\end{proof}

Referring to \cite[Section 1.6]{Koe2010}, let $\{h_i(x)\}_{i=0}^d$ denote the dual Hahn polynomials with parameters 
$$
(N,\gamma,\delta)=(d,-s-d-1,-r-d-1).
$$ 
Let ${\rm aff}_1:\R\to \R$ denote the affine transformation given by 
$x \mapsto  x+d(d+r+s+1)$
for all $x\in \R$. 
Then 
$$
h_i(x)=u_i({\rm aff}_1(x)) 
\qquad 
(0\leq i\leq d).
$$
Therefore $\{u_i(x)\}_{i=0}^d$ are essentially the dual Hahn polynomials.
By \cite[Section 1.6]{Askeyscheme} the dual Hahn polynomials $\{u_i(x)\}_{i=0}^d$ are a family of discrete orthogonal polynomials on the points $\{\theta_i\}_{i=0}^d$ with weights $\{k^*_i\}_{i=0}^d$. More precisely $\{u_i(x)\}_{i=0}^d$ obey the following orthogonality relation:
\begin{gather}
\label{orth_dualHahn}
\sum_{h=0}^d u_i(\theta_h) u_j(\theta_h) k^*_h=
\left\{
\begin{array}{ll}
0 
\qquad 
&\hbox{if $i\not=j$},
\\
\nu/k_i
\qquad 
&\hbox{if $i=j$}
\end{array}
\right.
\qquad 
(0\leq i,j\leq d).
\end{gather}
For convenience we always set each of  $b_d,b_d^*,c_0,c_0^*$ as the zero scalar and let $\theta_{-1}, \theta_{d+1}$ be any scalars and $u_{-1}(x), u_{d+1}(x)$ be any polynomials. 
The three-term recurrence relation for the dual Hahn polynomials $\{u_i(x)\}_{i=0}^d$ is 
\begin{gather}
\label{3term_dualHahn}
x u_i(x)=b_i u_{i+1}(x)+a_i u_i(x)+c_i u_{i-1}(x)
\qquad 
(0\leq i\leq d).
\end{gather}
The difference equations for the dual Hahn polynomials $\{u_i(x)\}_{i=0}^d$ are
\begin{gather}
\label{diff_dualHahn}
\theta^*_i u_i(\theta_j)=
b^*_ju_i(\theta_{j+1})
+a^*_j u_i(\theta_j)
+c^*_j u_i(\theta_{j-1})
\qquad 
(0\leq i,j\leq d).
\end{gather}

Next let us recall the notion of Leonard pairs. A square matrix is said to be {\it tridiagonal} whenever each nonzero entry lies on either the diagonal,
the subdiagonal, or the superdiagonal. A tridiagonal square matrix is said to be {\it irreducible}
whenever each entry on the subdiagonal is nonzero and each entry on the superdiagonal is
nonzero.

\begin{defn}
[Definition 1.1, \cite{Askeyscheme}]
\label{defn:LP}
Let $V$ denote a nonzero finite-dimensional vector space. A {\it Leonard pair} on $V$ consists of an ordered pair of linear transformations $A:V\to V$ and $A^*:V\to V$ that satisfy the following conditions:
\begin{enumerate}
\item There exists an ordered basis for $V$ with respect to which the matrix representing $A$ is diagonal and the matrix representing $A^*$ is irreducible tridiagonal.

\item There exists an ordered basis for $V$ with respect to which the matrix representing $A^*$ is diagonal and the matrix representing $A$ is irreducible tridiagonal.
\end{enumerate}
In this case, all eigenvalues of $A$ and $A^*$ are simple.
\end{defn}

Applying the three-term recurrence relation and difference equations, every family of discrete orthogonal polynomials in the Askey scheme has its own Leonard pair. The Leonard pair $L, L^*$ corresponding to the dual Hahn polynomial $\{u_i(x)\}_{i=0}^d$ is as follows:  
Throughout this paper, let $L$ denote the linear transformation on $\Pd$ given by  
\begin{eqnarray*}
(Lf)(\theta_i)  &= & 
\theta_i f(\theta_i)
\qquad 
(0\leq i\leq d)
\end{eqnarray*}
for all $f\in \Pd$.
Throughout this paper, let $L^*$ denote the linear transformation on $\Pd$ given by   
\begin{eqnarray*}
(L^*f)(\theta_i)  &= & 
b_i^* f(\theta_{i+1})
+a_i^* f(\theta_i)
+c_i^* f(\theta_{i-1})
\qquad 
(0\leq i\leq d)
\end{eqnarray*}
for all $f\in \Pd$. Applying (\ref{3term_dualHahn}) and  (\ref{diff_dualHahn}) yields that the matrices representing $L$ and $L^*$ with respect to the ordered basis $\{u_i\}_{i=0}^d$ for $\Pd$ are 
\begin{gather}
\label{L_tri&L*_diag}
\begin{pmatrix}
a_0 &c_1 & & &{\bf 0}
\\
b_0 &a_1 &c_2
\\
 &b_1 &a_2 &\ddots
\\
 & &\ddots  &\ddots  &c_d
\\
{\bf 0} & & &b_{d-1} &a_d
\end{pmatrix},
\qquad
\quad
\begin{pmatrix}
\theta_0^* & & & &{\bf 0}
\\
 &\theta_1^* &
\\
 & &\theta_2^* &
\\
 & &  &\ddots  &
\\
{\bf 0} & & & &\theta_d^*
\end{pmatrix},
\end{gather}
respectively. 
By (\ref{bi}), (\ref{ci}) and since $r,s\in (-1,\infty)$ the scalars $\{b_i\}_{i=0}^{d-1}$ and $\{c_i\}_{i=1}^d$ are nonzero. Hence the first square matrix in (\ref{L_tri&L*_diag}) is irreducible tridiagonal.
Since the scalars $\{k_i^*\}_{i=0}^d$ are nonzero, 
there are uniquely polynomials $\{u_i^*(x)\}_{i=0}^d$ with $u_i^*\in \Pd$ for all $i=0,1,\ldots,d$ satisfying 
\begin{gather*}
u_i^*(\theta_j)
=
\left\{
\begin{array}{ll}
0 
\qquad &\hbox{if $j\not=i$},
\\
\nu/k_i^*
\qquad &\hbox{if $j=i$}
\end{array}
\right.
\qquad 
(0\leq i,j\leq d).
\end{gather*}
Since the scalar $\nu$ is nonzero, the elements $\{u_i^*\}_{i=0}^d$ form a basis for $\Pd$. 
Using (\ref{ki*-1}) yields that the matrices representing $L$ and $L^*$ with respect to the ordered basis $\{u_i^*\}_{i=0}^d$ for $\Pd$ are 
\begin{gather}
\label{L_diag&L*_tri}
\begin{pmatrix}
\theta_0 & & & &{\bf 0}
\\
 &\theta_1 &
\\
 & &\theta_2 &
\\
 & &  &\ddots  &
\\
{\bf 0} & & & &\theta_d
\end{pmatrix},
\qquad 
\quad
\begin{pmatrix}
a_0^* &c_1^* & & &{\bf 0}
\\
b_0^* &a_1^* &c_2^*
\\
 &b_1^* &a_2^* &\ddots
\\
 & &\ddots  &\ddots  &c_d^*
\\
{\bf 0} & & &b_{d-1}^* &a_d^*
\end{pmatrix},
\end{gather}
respectively. 
By (\ref{bi*}), (\ref{ci*}) and since $r,s\in (-1,\infty)$ the scalars $\{b_i^*\}_{i=0}^{d-1}$ and $\{c_i^*\}_{i=1}^d$ are nonzero. Hence the second square matrix in (\ref{L_diag&L*_tri}) is irreducible tridiagonal.
Therefore $L,L^*$ is indeed a Leonard pair. 
The first main result of this paper is as follows:

\begin{thm}
\label{thm:LP}
Suppose that $r\not=0$ and $r+s=0$. Then $L,(L^*+\frac{r-d}{2})^2$ is a Leonard pair.
\end{thm}

Recall the Racah polynomials from \cite[Section 1.2]{Koe2010}. Applying the theory of Leonard pairs to Theorem \ref{thm:LP} yields the following result:

\begin{thm}
\label{thm:Racah}
Suppose that $r\not=0$ and $r+s=0$. 
Then $\{u_i(x)\}_{i=0}^d$ are not only the dual Hahn polynomials but also the Racah polynomials on the same points $\{\theta_i\}_{i=0}^d$ with the same weights $\{k_i^*\}_{i=0}^d$.
\end{thm}

As everyone knows, the dual Hahn polynomials are limiting cases of the Racah polynomials \cite[Section 2.2]{Koe2010}. Theorem \ref{thm:Racah} reveals another relationship between the dual Hahn polynomials and the Racah polynomials.
The hypothesis of Theorems \ref{thm:LP} and \ref{thm:Racah} along with the requisite $r,s\in (-1,\infty)$ is equivalent to 
$$
r\in (-1,0)\cup (0,1)
\quad 
\hbox{and}
\quad 
s=-r.
$$
When $r,s\in (-1,\infty)$ is replaced by $r,s\in (-\infty,-d)$,   	the points $\{\theta_i\}_{i=0}^d$ are still distinct and the weights $\{k_i^*\}_{i=0}^d$ are positive; however there are no  scalars $r$ and $s$ satisfying $r+s=0$.

The proofs for Theorems \ref{thm:LP} and \ref{thm:Racah} are given in the following two sections. 
In the final section we explain the relationship among $\U$, $\frac{1}{2}H(D,2)$ and Theorem \ref{thm:LP}.

\section{Proof for Theorem \ref{thm:LP}}\label{sec:LP}

For any linear transformation $T$ on $\Pd$, let $[T]$ denote the matrix representing $T$ with respect to the ordered basis $\{u_i^*\}_{i=0}^d$ for $\Pd$. 
Following the traditional convention in the field, the rows and columns of a $(d+1)\times (d+1)$ matrix are indexed by $0,1,\ldots,d$. 
Throughout this section the notation $\lambda$ stands for a real scalar.

\begin{lem}
\label{lem1:L*2}
For any integers $i,j$ with $0\leq i,j\leq d$, the $(i,j)$-entry of $[(L^*+\lambda)^2]$ is equal to 
\begin{gather*}
\left\{
\begin{array}{ll}
c_{i+1}^* c_{i+2}^*
\qquad &\hbox{if $j-i=2$},
\\
c^*_{i+1}(2\lambda+a_i^*+a^*_{i+1})
\qquad &\hbox{if $j-i=1$},
\\
b^*_{i} c^*_{i+1} 
+b^*_{i-1} c^*_{i}
+(\lambda+a^*_{i})^2
\qquad &\hbox{if $i=j$},
\\
b^*_{i-1}(2\lambda+a^*_i+a^*_{i-1})
\qquad &\hbox{if $i-j=1$},
\\
b_{i-1}^* b_{i-2}^*
\qquad &\hbox{if $i-j=2$},
\\
0
\qquad &\hbox{else},
\end{array}
\right.
\end{gather*}
where $b_{-1}^*$ and $c_{d+1}^*$ are interpreted as any scalars.
\end{lem}
\begin{proof}
It is routine to verify the lemma by using (\ref{L_diag&L*_tri}).
\end{proof}

Recall that 
the scalars $\{b_i^*\}_{i=0}^{d-1}$ and $\{c_i^*\}_{i=1}^d$ are nonzero. 

\begin{lem}
\label{lem2:L*2}
For any integer $i$ with $0\leq i\leq d-2$ the following statements are true:
\begin{enumerate}
\item $[(L^*+\lambda)^2]_{i,i+2}\not=0$.

\item $[(L^*+\lambda)^2]_{i+2,i}\not=0$.
\end{enumerate}
\end{lem}
\begin{proof}
Immediate from Lemma \ref{lem1:L*2}.
\end{proof}

\begin{lem}
\label{lem3:L*2}
For any integer $i$ with $0\leq i\leq d-1$ the following conditions are equivalent: 
\begin{enumerate}
\item $[(L^*+\lambda)^2]_{i,i+1}=0$.

\item $[(L^*+\lambda)^2]_{i+1,i}=0$.

\item $2\lambda+a^*_i+a^*_{i+1}=0$.
\end{enumerate}
\end{lem}
\begin{proof}
Immediate from Lemma \ref{lem1:L*2}.
\end{proof}

\begin{lem}
\label{lem:colsum}
Each column sum of $[(L^*+\lambda)^2]$ is equal to 
$\lambda^2$.
\end{lem}
\begin{proof}
Using (\ref{ai*}) and Lemma \ref{lem1:L*2} it is routine to verify the lemma.
\end{proof}

\begin{prop}
\label{prop:L&L*2}
Assume that $d\geq 1$.
Let $\{\bar{u}_i^*\}_{i=0}^d$ denote an ordering of $\{u_i^*\}_{i=0}^d$. 
Then the matrix representing $(L^*+\lambda)^2$ with respect to the ordered basis $\{\bar{u}_i^*\}_{i=0}^d$ for $\Pd$ is irreducible tridiagonal if and only if one of the following conditions {\rm (i)} and {\rm (ii)} holds:
\begin{enumerate}
\item $2\lambda+a_{d-1}^*+a_{d}^*\not=0$ and  $2\lambda+a_i^*+a_{i+1}^*=0$ for all $i=0,1,\ldots,d-2$. 
In addition 
\begin{align}
\label{baru:1}
\bar{u}_i^*
&=\left\{
\begin{array}{ll}
u_{2i}^*
\quad 
&\hbox{if $i=0,1,\ldots,\lfloor \frac{d}{2}\rfloor$},
\\
u_{2(d-i)+1}^*
\quad 
&\hbox{if $i=\lfloor\frac{d}{2}\rfloor+1,\lfloor\frac{d}{2}\rfloor+2,\ldots,d$}
\end{array}
\right.
\qquad 
(0\leq i\leq d)
\end{align}
or 
\begin{align}
\label{baru:2}
\bar{u}_i^*
&=\left\{
\begin{array}{ll}
u_{2i+1}^*
\quad 
&\hbox{if $i=0,1,\ldots,\lceil\frac{d}{2}\rceil-1$},
\\
u_{2(d-i)}^*
\quad 
&\hbox{if $i=\lceil \frac{d}{2}\rceil,\lceil \frac{d}{2}\rceil+1,\ldots,d$}
\end{array}
\right.
\qquad 
(0\leq i\leq d).
\end{align}

\item 
$2\lambda+a_0^*+a_1^*\not=0$ and  $2\lambda+a_i^*+a_{i+1}^*=0$ for all $i=1,2\ldots,d-1$. 
In addition  
\begin{align}
\label{baru:3}
\bar{u}_i^*
&=\left\{
\begin{array}{ll}
u_{d-2i}^*
\quad 
&\hbox{if $i=0,1,\ldots,\lfloor \frac{d}{2}\rfloor$},
\\
u_{2i-d-1}^*
\quad 
&\hbox{if $i=\lfloor\frac{d}{2}\rfloor+1,\lfloor\frac{d}{2}\rfloor+2,\ldots,d$}
\end{array}
\right.
\qquad 
(0\leq i\leq d)
\end{align}
or 
\begin{align}
\label{baru:4}
\bar{u}_i^*
&=\left\{
\begin{array}{ll}
u_{d-2i-1}^*
\quad 
&\hbox{if $i=0,1,\ldots,\lceil\frac{d}{2}\rceil-1$},
\\
u_{2i-d}^*
\quad 
&\hbox{if $i=\lceil \frac{d}{2}\rceil,\lceil \frac{d}{2}\rceil+1,\ldots,d$}
\end{array}
\right.
\qquad 
(0\leq i\leq d).
\end{align}
\end{enumerate}
\end{prop}
\begin{proof}
($\Leftarrow$): It is routine to verify the ``if'' part by using Lemmas \ref{lem1:L*2}--\ref{lem3:L*2}.

($\Rightarrow$): 
Suppose that the matrix representing $(L^*+\lambda)^2$ with respect to the ordered basis $\{\bar{u}_i^*\}_{i=0}^d$ for $\Pd$ is irreducible tridiagonal. 
By Lemma \ref{lem1:L*2}  any two consecutive polynomials $u^*_k$ and $u^*_{\ell}$ in the sequence $\{\bar{u}_i^*\}_{i=0}^d$ satisfy $1\leq |k-\ell|\leq 2$. 
By Lemma \ref{lem2:L*2} 
the sequence $\{\bar{u}^*_i\}_{i=0}^d$ contains the contiguous subsequence
\begin{gather*}
u^*_0, u^*_2, u^*_4, \ldots, u^*_{2\lfloor \frac{d}{2}\rfloor}
\end{gather*}
or its reverse sequence. Also $\{\bar{u}^*_i\}_{i=0}^d$ contains the contiguous subsequence
\begin{gather*}
u^*_1, u^*_3, u^*_5, \ldots, u^*_{2\lceil \frac{d}{2}\rceil-1}
\end{gather*}
or its reverse sequence. 
Based on the above comments the sequence $\{\bar{u}^*_i\}_{i=0}^d$ must be one of (\ref{baru:1})--(\ref{baru:4}). Combined with Lemma \ref{lem3:L*2} the ``only if'' part follows.
\end{proof}

\begin{lem}
\label{lem:a*i+a*i+1}
The scalar $a^*_i+a^*_{i+1}$ is equal to 
\begin{gather*}
\left\{
\begin{array}{ll}
\displaystyle 
d-\frac{r-s}{2}+\frac{(r-s)(r+s)(2d+r+s+2)}{2(2d-2i+r+s-2)(2d-2i+r+s+2)}
\qquad &\hbox{if $i=0,1,\ldots,d-2$},
\\
\displaystyle
d+\frac{(d-1)(r-s)}{r+s+4}
\qquad &\hbox{if $i=d-1$.}
\end{array}
\right.
\end{gather*}
\end{lem}
\begin{proof}
Evaluate $a^*_i+a^*_{i+1}$ by applying (\ref{thetai*}) and  (\ref{bi*})--(\ref{ai*}). 
\end{proof}

\begin{lem}
\label{lem3:a*i}
\begin{enumerate}
\item Suppose that $r+s=0$. Then 
$$
a^*_i=\frac{d-r}{2} \qquad (0\leq i\leq d-1), 
\qquad
a^*_d=\frac{d(r+1)}{2}. 
$$

\item Suppose that $r-s=0$. Then 
$$
a^*_i=\frac{d}{2}
\qquad 
(0\leq i\leq d).
$$
\end{enumerate}
\end{lem}
\begin{proof}
(i): By (\ref{bi*}) and (\ref{ci*}) the scalars 
\begin{align*}
b^*_i&=\frac{(2d-i+1)(i-d+r)}{2(2d-2i+1)}
\qquad (0\leq i\leq d-1),
\\
c^*_i&=\frac{i(i-d-r-1)}{2(2d-2i+1)}
\qquad (1\leq i\leq d)
\end{align*}
under the hypothesis $r+s=0$. Combined with (\ref{ai*}) the statement (i) follows.

(ii): By (\ref{bi*}) and (\ref{ci*}) the scalars 
\begin{align*}
b^*_i&=-\frac{(d-i)(2d-i+2r+1)}{2(2d-2i+2r+1)}
\qquad (0\leq i\leq d-1),
\\
c^*_i&=-\frac{i(d-i+2r+1)}{2(2d-2i+2r+1)}
\qquad (1\leq i\leq d)
\end{align*}
under the hypothesis $r-s=0$. Combined with (\ref{ai*}) the statement (ii) follows.
\end{proof}

\begin{lem}
\label{lem1:a*i}
Suppose that $d\geq 2$. Then the following conditions are equivalent: 
\begin{enumerate}
\item $a^*_{d-2}+a^*_{d-1}=a^*_{d-1}+a^*_d$.

\item $r-s=0$.
\end{enumerate}
\end{lem}
\begin{proof}
(ii) $\Rightarrow$ (i): Immediate from Lemma \ref{lem3:a*i}(ii).

(i) $\Rightarrow$ (ii): Using Lemma \ref{lem:a*i+a*i+1} yields that $(a^*_{d-1}+a^*_d)-(a^*_{d-2}+a^*_{d-1})$ is equal to $r-s$ times 
\begin{gather}
\label{a*d-a*d-2}
\frac{(r+s+3)(2d+r+s+2)}{4(\frac{r+s}{2}+1)_3}.
\end{gather}
Since $r,s\in (-1,\infty)$ the fraction (\ref{a*d-a*d-2}) is positive. Therefore (i) implies (ii). 
\end{proof}

\begin{lem}
\label{lem2:a*i}
Suppose that $d\geq 3$. Then the following conditions are equivalent:
\begin{enumerate}
\item $a^*_i+a^*_{i+1}=a^*_{i+1}+a^*_{i+2}$ for each integer $i$ with $0\leq i\leq d-3$.

\item $a^*_i+a^*_{i+1}=a^*_{i+1}+a^*_{i+2}$ for some integer $i$ with $0\leq i\leq d-3$.

\item $r-s=0$ or $r+s=0$.
\end{enumerate}
\end{lem}
\begin{proof}
Let $i$ be an integer with $0\leq i\leq d-3$. Using Lemma \ref{lem:a*i+a*i+1} yields that $(a^*_{i+2}+a^*_{i+1})-(a^*_{i+1}+a^*_i)$ is equal to $(r-s)(r+s)$ times  
\begin{gather}
\label{a*i+2-a*i}
\frac{(2d+r+s+2)(2d-2i+r+s-1)}
{8(d-i+\frac{r+s}{2}-2)_4}.
\end{gather}
Since $r,s\in (-1,\infty)$ the fraction (\ref{a*i+2-a*i}) is positive. Therefore (i)--(iii) are equivalent. 
\end{proof}

Now we can give the conditions for $L, (L^*+\lambda)^2$ as a Leonard pair. 
By (\ref{L_tri&L*_diag}) the matrices representing $L$ and $(L^*+\lambda)^2$ with respect to the ordered basis $\{u_i\}_{i=0}^d$ for $\Pd$ are irreducible tridiagonal and diagonal, respectively. Thus Definition \ref{defn:LP}(ii) holds for the pair $L, (L^*+\lambda)^2$.

\begin{cor}
Suppose that $d=1$. Then the following conditions are equivalent:
\begin{enumerate}
\item $L,(L^*+\lambda)^2$ is a Leonard pair.

\item $2\lambda\not=-1$.
\end{enumerate}
\end{cor}
\begin{proof}
By Lemma \ref{lem:a*i+a*i+1} the scalar  $a_0^*+a_1^*=1$. Combined with Proposition \ref{prop:L&L*2} this yields the corollary.
\end{proof}

\begin{cor}
Suppose that $d=2$. Then the following conditions are equivalent:
\begin{enumerate}
\item $L,(L^*+\lambda)^2$ is a Leonard pair.

\item $r\not=s$ and $
2(\lambda+1)\in\{\frac{r-s}{r+s+2},\frac{s-r}{r+s+4}\}$.
\end{enumerate}
\end{cor}
\begin{proof}
By Proposition \ref{prop:L&L*2} the statement (i) is true if and only if 
\begin{gather}
a_0^*+a_1^*\not=a_1^*+a_2^*,
\label{d=2_cond1}
\\
-2\lambda\in \{a^*_0+a^*_1,a^*_1+a^*_2\}.
\label{d=2_cond2}
\end{gather} 
By Lemma \ref{lem1:a*i}  the condition (\ref{d=2_cond1}) is equivalent to $r\not=s$. 
By Lemma \ref{lem:a*i+a*i+1} the scalars $a_0^*+a_1^*$ and $a_1^*+a_2^*$ are equal to $2-\frac{r-s}{r+s+2}$ and $2+\frac{r-s}{r+s+4}$, respectively. 
The condition (\ref{d=2_cond2}) is equivalent to $
2(\lambda+1)\in\{\frac{r-s}{r+s+2},\frac{s-r}{r+s+4}\}$. The corollary follows.
\end{proof}

\begin{thm}
\label{thm:LP&r,s,lambda}
Suppose that the following conditions hold:
\begin{enumerate}
\item $r\not=0$.

\item $r+s=0$.

\item $2\lambda=r-d$. 
\end{enumerate}
Then $L,(L^*+\lambda)^2$ is a Leonard pair. 
The converse is true provided that $d\geq 3$. 
\end{thm}
\begin{proof}
If $d=0$ there is nothing to prove. Thus we assume that $d\geq 1$. 
Suppose that (i)--(iii) hold. It follows from Lemma \ref{lem3:a*i} that 
\begin{align*}
2\lambda+a^*_{d-1}+a^*_d&=\frac{r(d+1)}{2}\not=0,
\\
2\lambda+a^*_i+a^*_{i+1}&=0
\qquad 
(0\leq i\leq d-2).
\end{align*} 
Combined with Proposition \ref{prop:L&L*2}, Definition \ref{defn:LP}(i) holds for the pair $L, (L^*+\lambda)^2$. The first assertion follows.

Suppose that $d\geq 3$ and $L,(L^*+\lambda)^2$ is a Leonard pair. By Definition \ref{defn:LP}(i) there exists an ordering $\{\bar{u}^*_i\}_{i=0}^d$ of the basis $\{u^*_i\}_{i=0}^d$ for $\Pd$ with respect to which the matrix representing $(L^*+\lambda)^2$ is irreducible tridiagonal. Thus either Proposition \ref{prop:L&L*2}(i) or Proposition \ref{prop:L&L*2}(ii) holds. Suppose that Proposition \ref{prop:L&L*2}(ii) holds. Then $r-s=0$ by Lemma \ref{lem1:a*i}. Combined with Lemma \ref{lem2:a*i} this implies that
$$
a_i^*+a_{i+1}^*=a_{i+1}^*+a_{i+2}^*
\qquad (0\leq i\leq d-2),
$$
a contradiction to Proposition \ref{prop:L&L*2}(ii).
Therefore Proposition \ref{prop:L&L*2}(i) holds. Combined with Lemmas \ref{lem3:a*i}--\ref{lem2:a*i} the conditions (i)--(iii) follow. The second assertion follows.
\end{proof}

\begin{proof}[Proof of Theorem \ref{thm:LP}] Immediate from Theorem \ref{thm:LP&r,s,lambda}. 
\end{proof}

Inspired by Theorem \ref{thm:LP} one might be interested in Leonard pairs with the following property: 

\begin{defn}
A Leonard pair $A,A^*$ is said to be {\it square-preserving}
whenever one of the following conditions holds:
\begin{enumerate}
\item $A^2,A^*$ is a Leonard pair. 

\item $A,{A^*}^2$ is a Leonard pair.
\end{enumerate}
\end{defn}

Recall the isomorphism of Leonard pairs.

\begin{defn}
[Definition 3.4, \cite{Askeyscheme}]
\label{defn:LPiso}
Assume that $V$ and $W$ are nonzero finite-dimensional vector spaces over the same field.
Let $A,A^*$ denote a Leonard pair on $V$. Let $B,B^*$ denote a Leonard pair on $W$. Then the Leonard pair $A,A^*$ is said to be {\it isomorphic} to the Leonard pair $B,B^*$ if there exists an invertible linear transformation $T:V\to W$ such that $T\circ A\circ T^{-1}=B$ and $T\circ A^*\circ T^{-1}=B^*$.
\end{defn}

When a Leonard pair $A,A^*$ is isomorphic to a Leonard pair $B,B^*$, the Leonard pair $A,A^*$ is square-preserving if and only if the Leonard pair $B,B^*$ is square-preserving. While sending the manuscript to Kazumasa Nomura, he proposed a problem like:

\begin{prob}
\label{prob:Kazumasa}
Please classify the square-preserving Leonard pairs up to isomorphism.
\end{prob}

\section{Proof for Theorem \ref{thm:Racah}}\label{sec:Racah}

Throughout this section we always assume that $r\not=0$ and $r+s=0$. In addition we will use the following parameters: 
\begin{align}
\bar{\theta}_i&=(d-2i)(d-2i+1)
\qquad (0\leq i\leq d),
\label{barthetai}
\\
\bar{\theta}_i^*&=
\textstyle 
(i+\frac{r-d}{2})^2
\qquad (0\leq i\leq d),
\label{bartheta*i}
\\
\bar{b}_i&=(d-i)(d-i-r)
\qquad (0\leq i\leq d-1),
\label{barbi}
\\
\bar{c}_i&=i(i+r) \qquad (1\leq i\leq d),
\label{barci}
\\
\bar{a}_i&=
\left\{
\begin{array}{ll}
\bar{\theta}_0-\bar{b}_0
\qquad &\hbox{if $i=0$},
\\
\bar{\theta}_0-\bar{b}_i-\bar{c}_i
\qquad &\hbox{if $i=1,2,\ldots,d-1$},
\\
\bar{\theta}_0-\bar{c}_d
\qquad &\hbox{if $i=d$}
\end{array}
\right.
\qquad (0\leq i\leq d),
\label{barai}
\\
\bar{k}_i&=
\frac{\bar{b}_0\bar{b}_1\cdots \bar{b}_{i-1}}
{\bar{c}_1\bar{c}_2\cdots \bar{c}_i}
\qquad (0\leq i\leq d),
\label{barki}
\\
\bar{\nu} &=\frac{(\bar{\theta}_0-\bar{\theta}_1)(\bar{\theta}_0-\bar{\theta}_2)\cdots (\bar{\theta}_0-\bar{\theta}_d)}{\bar{c}_1\bar{c}_2\cdots\bar{c}_d},
\label{barnu}
\\
\bar{b}_i^* &=
\frac{(d-i)(2d-2i+1)(d-2i-r-1)(d-2i-r)}{2(2d-4i-1)(2d-4i+1)}
\qquad (0\leq i\leq d-1),
\label{barb*i}
\\
\bar{c}_i^* &=
\frac{i(2i-1)(d-2i+r+1)(d-2i+r+2)}
{2(2d-4i+1)(2d-4i+3)}
\qquad (1\leq i\leq d),
\label{barc*i}
\\
\bar{a}_i^*&=\left\{
\begin{array}{ll}
\bar{\theta}_0^*-\bar{b}_0^*
\qquad &\hbox{if $i=0$},
\\
\bar{\theta}_0^*-\bar{b}_i^*-\bar{c}_i^*
\qquad &\hbox{if $i=1,2,\ldots,d-1$},
\\
\bar{\theta}_0^*-\bar{c}_d^*
\qquad &\hbox{if $i=d$}
\end{array}
\right.
\qquad (0\leq i\leq d),
\label{bara*i}
\\
\bar{k}_i^* &=
\frac{\bar{b}^*_0\bar{b}^*_1\cdots \bar{b}^*_{i-1}}
{\bar{c}^*_1\bar{c}^*_2\cdots \bar{c}^*_i}
\qquad (0\leq i\leq d).
\label{bark*i}
\end{align}

\begin{lem}
\label{lem:bartheta}
The following equations hold:
\begin{enumerate}
\item 
$
\bar{\theta}_i
=\left\{
\begin{array}{ll}
\theta_{2i}
\quad 
&\hbox{if $i=0,1,\ldots,\lfloor \frac{d}{2}\rfloor$},
\\
\theta_{2(d-i)+1}
\quad 
&\hbox{if $i=\lfloor\frac{d}{2}\rfloor+1,\lfloor\frac{d}{2}\rfloor+2,\ldots,d$}
\end{array}
\right.
\qquad 
(0\leq i\leq d).
$

\item 
$
\bar{\theta}^*_i=(\theta^*_i+\frac{r-d}{2})^2$ for all $i=0,1,\ldots,d$.
\end{enumerate}
\end{lem}
\begin{proof}
It is straightforward to verify (i) by using (\ref{thetai}) and (\ref{barthetai}). 
It is straightforward to verify (ii) by using (\ref{thetai*}) and  (\ref{bartheta*i}).
\end{proof}

\begin{lem}
\label{lem:barbck}
The following equations hold:
\begin{enumerate}
\item $\bar{b}_i=b_i$ for all $i=0,1,\ldots,d-1$.

\item $\bar{c}_i=c_i$ for all $i=1,2,\ldots,d$.

\item $\bar{a}_i=a_i$ for all $i=0,1,\ldots,d$.

\item $\bar{k}_i=k_i$ for all $i=0,1,\ldots,d$.

\item $\bar{\nu}=\nu$.
\end{enumerate}
\end{lem}
\begin{proof}
(i): Immediate from (\ref{bi}) and (\ref{barbi}).

(ii): Immediate from (\ref{ci}) and (\ref{barci}).

(iii): Immediate from (i), (ii), (\ref{ai}), (\ref{barai}) and Lemma \ref{lem:bartheta}(i).

(iv): Immediate from (i), (ii) and (\ref{ki-1}), (\ref{barki}).

(v):  Immediate from (ii), (\ref{nu-1}), (\ref{barnu}) and Lemma \ref{lem:bartheta}(i).
\end{proof}

\begin{lem}
\label{lem:barbck*}
The following equations hold:
\begin{enumerate}
\item
$\bar{b}^*_i
=\left\{
\begin{array}{ll}
b^*_{2i} b^*_{2i+1} 
\quad 
&\hbox{if $i=0,1,\ldots,\lfloor \frac{d}{2}\rfloor-1$},
\\
\frac{(d+1)r}{2}\cdot b^*_{d-1}
\quad 
&\hbox{if $i=\lfloor \frac{d}{2}\rfloor$ and $d$ is odd},
\\
\frac{(d+1)r}{2}\cdot c^*_d
\quad 
&\hbox{if $i=\lfloor \frac{d}{2}\rfloor$ and $d$ is even},
\\
c^*_{2(d-i)} c^*_{2(d-i)+1}
\quad 
&\hbox{if $i=\lfloor\frac{d}{2}\rfloor+1,\lfloor\frac{d}{2}\rfloor+2,\ldots,d-1$}
\end{array}
\right.
\qquad 
(0\leq i\leq d-1).
$

\item $
\bar{c}^*_i
=\left\{
\begin{array}{ll}
c^*_{2i} c^*_{2i-1}
\quad 
&\hbox{if $i=1,2,\ldots,\lfloor \frac{d}{2}\rfloor$},
\\
\frac{(d+1)r}{2}\cdot c^*_d
\quad 
&\hbox{if $i=\lfloor \frac{d}{2}\rfloor+1$ and $d$ is odd},
\\
\frac{(d+1)r}{2}\cdot b^*_{d-1}
\qquad
&\hbox{if $i=\lfloor \frac{d}{2}\rfloor+1$ and $d$ is even},
\\
b^*_{2(d-i+1)} b^*_{2(d-i)+1}
\quad 
&\hbox{if $i=\lfloor\frac{d}{2}\rfloor+2,\lfloor\frac{d}{2}\rfloor+2,\ldots,d$}
\end{array}
\right.
\qquad 
(1\leq i\leq d).
$

\item $\bar{k}^*_i=
\left\{
\begin{array}{ll}
k^*_{2i}
\quad 
&\hbox{if $i=0,1,\ldots,\lfloor \frac{d}{2}\rfloor$},
\\
k^*_{2(d-i)+1}
\quad 
&\hbox{if $i=\lfloor\frac{d}{2}\rfloor+1,\lfloor\frac{d}{2}\rfloor+2,\ldots,d$}
\end{array}
\right.
\qquad 
(0\leq i\leq d).
$
\end{enumerate}
\end{lem}
\begin{proof}
(i), (ii): It is routine to verify (i), (ii) by using (\ref{bi*}),  (\ref{ci*}) and (\ref{barb*i}), (\ref{barc*i}). 

(iii): Applying (i), (ii) to (\ref{bark*i}) yields that 
$$
\bar{k}^*_i
=
\left\{
\begin{array}{ll}
\displaystyle 
\frac{b^*_0 b^*_1\cdots b^*_{2i-1}}
{c^*_1 c^*_2\cdots c^*_{2i}}
\quad 
&\hbox{if $i=0,1,\ldots,\floor{\frac{d}{2}}$},
\\
\displaystyle 
\frac{b^*_0b^*_1\cdots b^*_{d-1}}
{c^*_1c^*_2\cdots c^*_{d}}
\quad 
&\hbox{if $i=\lfloor \frac{d}{2}\rfloor+1$ and $d$ is odd},
\\
\displaystyle 
\frac{b^*_0b^*_1\cdots b^*_{d-2}}
{c^*_1c^*_2\cdots c^*_{d-1}}
\quad 
&\hbox{if $i=\lfloor \frac{d}{2}\rfloor+1$ and $d$ is even},
\\
\displaystyle 
\frac{b^*_0b^*_1\cdots b^*_{2(d-i)}}
{c^*_1c^*_2\cdots c^*_{2(d-i)+1}}
\quad 
&\hbox{if $i=\lfloor \frac{d}{2}\rfloor+2,\lfloor \frac{d}{2}\rfloor+3,\ldots,d$}
\end{array}
\right.
\qquad 
(0\leq i\leq d).
$$
The equation (iii) follows by applying (\ref{ki*-1}) to the right-hand side of the above equation. 
\end{proof}

Recall the matrix representations (\ref{L_tri&L*_diag}) of $L$ and $L^*$ with respect to the ordered basis $\{u_i\}_{i=0}^d$ for $\Pd$.
By Lemmas \ref{lem:bartheta}(ii) and \ref{lem:barbck}(i)--(iii) the matrices representing $L$ and $(L^*+\frac{r-d}{2})^2$ with respect to the ordered basis $\{u_i\}_{i=0}^d$ for $\Pd$ are 
\begin{gather}
\label{L_tri&barL*_diag}
\begin{pmatrix}
\bar{a}_0 &\bar{c}_1 & & &{\bf 0}
\\
\bar{b}_0 &\bar{a}_1 &\bar{c}_2
\\
 &\bar{b}_1 &\bar{a}_2 &\ddots
\\
 & &\ddots  &\ddots  &\bar{c}_d
\\
{\bf 0} & & &\bar{b}_{d-1} &\bar{a}_d
\end{pmatrix},
\qquad
\quad
\begin{pmatrix}
\bar{\theta}_0^* & & & &{\bf 0}
\\
 &\bar{\theta}_1^* &
\\
 & &\bar{\theta}_2^* &
\\
 & &  &\ddots  &
\\
{\bf 0} & & & &\bar{\theta}_d^*
\end{pmatrix},
\end{gather}
respectively. 
Recall the matrix representations (\ref{L_diag&L*_tri})
of $L$ and $L^*$ with respect to the ordered basis $\{u^*_i\}_{i=0}^d$ for $\Pd$. 
Let $\{\bar{u}_i^*\}_{i=0}^d$ denote the ordering (\ref{baru:1}) of $\{u^*_i\}_{i=0}^d$. 
Applying Lemmas \ref{lem1:L*2}, \ref{lem:colsum}, \ref{lem:bartheta}(i) and \ref{lem:barbck*}(i), (ii) the matrices representing $L$ and $(L^*+\frac{r-d}{2})^2$ with respect to the ordered basis $\{\bar{u}_i^*\}_{i=0}^d$ for $\Pd$ are 
\begin{gather}
\label{L_diag&barL*_tri}
\begin{pmatrix}
\bar{\theta}_0 & & & &{\bf 0}
\\
 &\bar{\theta}_1 &
\\
 & &\bar{\theta}_2 &
\\
 & &  &\ddots  &
\\
{\bf 0} & & & &\bar{\theta}_d
\end{pmatrix},
\qquad 
\quad
\begin{pmatrix}
\bar{a}_0^* &\bar{c}_1^* & & &{\bf 0}
\\
\bar{b}_0^* &\bar{a}_1^* &\bar{c}_2^*
\\
 &\bar{b}_1^* &\bar{a}_2^* &\ddots
\\
 & &\ddots  &\ddots  &\bar{c}_d^*
\\
{\bf 0} & & &\bar{b}_{d-1}^* &\bar{a}_d^*
\end{pmatrix},
\end{gather}
respectively.

By Lemma \ref{lem:barbck}(i)--(iii) 
the recurrence (\ref{3term_dualHahn}) can be written as follows:
$$
xu_i(x)=
\bar{b}_i u_{i+1}(x)
+
\bar{a}_i u_i(x)
+
\bar{c}_i u_{i-1}(x)
\qquad 
(0\leq i\leq d),
$$
where $\bar{b}_d$ and $\bar{c}_0$ are interpreted as the zero scalar. Combined with the matrix representations (\ref{L_tri&barL*_diag}) and (\ref{L_diag&barL*_tri}) of $L$ and $(L^*+\frac{r-d}{2})^2$, the polynomials $\{u_i(x)\}_{i=0}^d$ have the following alternative expressions:

\begin{lem}
\label{lem:4F3}
The polynomials $\{u_i(x)\}_{i=0}^d$ satisfy the following equations:
$$
u_i(\bar{\theta}_j) ={}_4F_3
\left(\genfrac..{0pt}{}{-i,
i-d+r,-j, 
j-d-\frac{1}{2}}
{-d,
\frac{r-d}{2},
\frac{r-d+1}{2}}
\,\Bigg|\,
1\right)
\qquad 
(0\leq i,j\leq d).
$$
\end{lem}
\begin{proof}
Let 
\begin{align*}
\bar{\varphi}_i
&=\bar{b}_{i-1}
\cdot 
\frac{
(\bar{\theta}^*_{i}-\bar{\theta}^*_0)
(\bar{\theta}^*_{i}-\bar{\theta}^*_1)
\cdots
(\bar{\theta}^*_{i}-\bar{\theta}^*_{i-1})}
{(\bar{\theta}^*_{i-1}-\bar{\theta}^*_0)
(\bar{\theta}^*_{i-1}-\bar{\theta}^*_1)
\cdots
(\bar{\theta}^*_{i-1}-\bar{\theta}^*_{i-2})}
\qquad 
(1\leq i\leq d).
\end{align*}
Substituting (\ref{barthetai}) and (\ref{barbi}) into the above equation yields that 
\begin{align}
\label{barvarphi}
\bar{\varphi}_i &=i(i-d-1)
(d-2i-r+1)(d-2i-r+2)
\qquad 
(1\leq i\leq d).
\end{align}
By \cite[Theorems 23.2 and 23.5]{Askeyscheme}, for any integers $i,j$ with $0\leq i,j\leq d$ the following equation holds:
$$
u_i(\bar{\theta}_j)
=\sum_{h=0}^i
\frac{
(\bar{\theta}^*_i-\bar{\theta}^*_0)
(\bar{\theta}^*_i-\bar{\theta}^*_1)
\cdots
(\bar{\theta}^*_i-\bar{\theta}^*_{h-1})
\cdot 
(\bar{\theta}_j-\bar{\theta}_0)
(\bar{\theta}_j-\bar{\theta}_1)
\cdots 
(\bar{\theta}_j-\bar{\theta}_{h-1})}
{\bar{\varphi}_1\bar{\varphi}_2\cdots\bar{\varphi}_h}.
$$
The lemma follows by substituting (\ref{barthetai}), (\ref{bartheta*i}), (\ref{barvarphi}) into the above equation and expressing the resulting equation in terms of ${}_4F_3$.
\end{proof}

\begin{proof}[Proof of Theorem \ref{thm:Racah}] 
Referring to \cite[Section 1.2]{Koe2010}, let $\{r_i(x)\}_{i=0}^d$ denote the Racah polynomials with parameters 
$$
(N,\alpha,\beta,\gamma,\delta)=
\textstyle
(d,-d-1,r,\frac{r-d-1}{2},-\frac{r+d}{2}-1)
$$
Let ${\rm aff}_2:\R\to \R$ denote the affine transformation given by 
$x \mapsto  4x+d(d+1)$ for all $x\in \R$. 
By Lemma \ref{lem:4F3} it is routine to verify that
$$
r_i(x)=u_i({\rm aff}_2(x))
\qquad 
(0\leq i\leq d).
$$
Therefore $\{u_i(x)\}_{i=0}^d$ are also the Racah polynomials.
By \cite[Theorem 18.3]{Askeyscheme} the orthogonality relation for the Racah polynomials $\{u_i(x)\}_{i=0}^d$ is as follows:
$$
\sum_{h=0}^d u_i(\bar{\theta}_h) u_j(\bar{\theta}_h)
\bar{k}_h^*
=
\left\{
\begin{array}{ll}
0 
\qquad 
&\hbox{if $i\not=j$},
\\
\bar{\nu}/\bar{k}_i
\qquad 
&\hbox{if $i=j$}
\end{array}
\right.
\qquad 
(0\leq i,j\leq d).
$$
By Lemmas \ref{lem:bartheta}(i), \ref{lem:barbck}(iv), (v) and \ref{lem:barbck*}(iii), the above orthogonality relation is identical to the orthogonality relation (\ref{orth_dualHahn}) for the dual Hahn polynomials $\{u_i(x)\}_{i=0}^d$. The result follows. 
\end{proof}


\section{The connection to $\U$ and $\frac{1}{2}H(D,2)$}\label{sec:LL*&U}

While showing an early version of the manuscript to Paul Terwilliger, he mentioned that when $D$ is odd, the distance-regular graph $\frac{1}{2}H(D,2)$ is almost dual-bipartite  \cite{nomura2025-2}, which is closely related to the so-called dual almost bipartite Leonard pairs \cite{Masuda2024,Brown2013}.

Let $A,A^*$ denote a Leonard pair on a nonzero finite-dimensional vector space $V$. The Leonard pair $A,A^*$ is called {\it almost bipartite} if there exists an ordered basis for $V$ with respect to which the matrix representing $A^*$ is diagonal and the matrix representing $A$ is irreducible diagonal and all diagonal entries except the last one are zero. The Leonard pair $A,A^*$ is called {\it dual almost bipartite} whenever the Leonard pair $A^*,A$ is almost bipartite. 
Note that every (dual) almost bipartite Leonard pair is square-preserving.

\begin{cor}
\label{cor:almost}
Suppose that $r\not=0$ and $r+s=0$. Then the Leonard pair $L,L^*+\frac{r-d}{2}$ is dual almost bipartite.
\end{cor}
\begin{proof}
Immediate from Lemma \ref{lem3:a*i}(i).
\end{proof}

Recall that $\U$ is a unital associative algebra over $\C$ generated by $E,F,H$ subject to the relations
\begin{gather*}
[H,E]=2E, 
\qquad
[H,F]=-2F,
\qquad
[E,F]=H. 
\end{gather*}
The Casimir element of $\U$ is defined as 
$$
\Lambda=EF+FE+\frac{H^2}{2}.
$$
Let $\U_e$ denote the unital subalgebra of $\U$ generated by $E^2,F^2,\Lambda,H$. The algebra $\U_e$ is called the {\it even subalgebra of $\U$} \cite[Definition 1.4 and Theorem 3.4]{halved:2023}.

\begin{lem}
\label{lem:Ln(0)}
For each integer $n\geq 0$ there exists a $(\lfloor\frac{n}{2}\rfloor+1)$-dimensional irreducible $\U_e$-module $L_n^{(0)}$ satisfying the following conditions:
\begin{enumerate}
\item There exists a basis $\{v_i^{(0)}\}_{i=0}^{\lfloor\frac{n}{2}\rfloor}$ for $L_n^{(0)}$ such that 
\begin{align*}
E^2 v_i^{(0)} &=2i(2i-1) v_{i-1}^{(0)}
\quad 
\textstyle (1\leq i\leq \lfloor\frac{n}{2}\rfloor),
\qquad 
E^2 v_0^{(0)} =0,
\\
F^2 v_i^{(0)} &= (n-2i)(n-2i-1) v_{i+1}^{(0)}
\quad 
\textstyle (0\leq i\leq \lfloor\frac{n}{2}\rfloor-1),
\qquad 
F^2 v_{\lfloor\frac{n}{2}\rfloor}^{(0)} =0,
\\
H v_i^{(0)} &=(n-4i) v_i^{(0)}
\quad 
\textstyle (0\leq i\leq \lfloor\frac{n}{2}\rfloor).
\end{align*}

\item The element $\Lambda$ acts on $L_n^{(0)}$ as scalar multiplication by $\frac{n(n+2)}{2}$.
\end{enumerate}
\end{lem}
\begin{proof}
The $\U_e$-module $L_n^{(0)}$ comes from \cite[Definition 5.2(i)]{halved:2023}. 
Let $\{u_i^{(0)}\}_{i=0}^{\lfloor \frac{n}{2}\rfloor}$ denote the basis for $L_n^{(0)}$ from \cite[Lemma 5.4(i)]{halved:2023}.
Then the vectors
$$
v_i^{(0)}={n \choose 2i}^{-1} u_i^{(0)}
\qquad \textstyle(0\leq i\leq \lfloor \frac{n}{2}\rfloor)
$$
give the basis for $L_n^{(0)}$ described in the statement (i). 
The statement (ii) is exactly \cite[Lemma 5.4(ii)]{halved:2023}. 
The irreducibility of the $\U_e$-module $L_n^{(0)}$ is shown in \cite[Lemma 5.5]{halved:2023}. 
\end{proof}

\begin{lem}
\label{lem:Ln(1)}
For each integer $n\geq 1$ there exists a $(\lfloor\frac{n+1}{2}\rfloor)$-dimensional irreducible $\U_e$-module $L_n^{(1)}$ satisfying the following conditions:
\begin{enumerate}
\item There exists a basis $\{v_i^{(1)}\}_{i=0}^{\lfloor\frac{n-1}{2}\rfloor}$ for $L_n^{(1)}$ such that 
\begin{align*}
E^2 v_i^{(1)} &=2i(2i+1) v_{i-1}^{(1)}
\quad 
\textstyle (1\leq i\leq \lfloor\frac{n-1}{2}\rfloor),
\qquad 
E^2 v_0^{(1)} =0,
\\
F^2 v_i^{(1)} &= (n-2i-1)(n-2i-2) v_{i+1}^{(1)}
\quad 
\textstyle (0\leq i\leq \lfloor\frac{n-3}{2}\rfloor),
\qquad 
F^2 v_{\lfloor\frac{n-1}{2}\rfloor}^{(1)} =0,
\\
H v_i^{(1)} &=(n-4i-2) v_i^{(1)}
\quad 
\textstyle (0\leq i\leq \lfloor\frac{n-1}{2}\rfloor).
\end{align*}

\item The element $\Lambda$ acts on $L_n^{(1)}$ as scalar multiplication by $\frac{n(n+2)}{2}$.
\end{enumerate}
\end{lem}
\begin{proof}
The $\U_e$-module $L_n^{(1)}$ comes from \cite[Definition 5.2(ii)]{halved:2023}. 
Let $\{u_i^{(1)}\}_{i=0}^{\lfloor \frac{n-1}{2}\rfloor}$ denote the basis for $L_n^{(1)}$ from \cite[Lemma 5.7(i)]{halved:2023}.
Then the vectors
$$
v_i^{(1)}=n\cdot {n \choose 2i+1}^{-1} u_i^{(1)}
\qquad \textstyle(0\leq i\leq \lfloor \frac{n-1}{2}\rfloor)
$$
give the basis for $L_n^{(1)}$ described in the statement (i). 
The statement (ii) is exactly \cite[Lemma 5.7(ii)]{halved:2023}. 
The irreducibility of the $\U_e$-module $L_n^{(1)}$ is shown in \cite[Lemma 5.8]{halved:2023}. 
\end{proof}

The $\U_e$-modules $L_n^{(0)}$ for all integers $n\geq 0$ and the $\U_e$-modules $L_n^{(1)}$ for all integers $n\geq 1$ are mutually non-isomorphic \cite[Theorem 5.10]{halved:2023}.

\begin{exam}
\label{exam:LP_Ln(0)}
Assume that $n\geq 1$ is an odd integer. Set the parameters 
$$
(r,s,d)=\left(-\frac{1}{2},\frac{1}{2},\frac{n-1}{2}\right).
$$
The settings of $r$ and $s$ fit the hypothesis of Theorem \ref{thm:LP}. 
Using Lemma \ref{lem:Ln(0)} yields that the matrices representing
\begin{gather}
\label{LP:Ln(0)}
\frac{E^2+F^2+\Lambda-1}{4}-\frac{H^2}{8}, \qquad \frac{n-H}{4}
\end{gather}
with respect to the ordered basis $\{v_i^{(0)}\}_{i=0}^d$ for $L_n^{(0)}$ are exactly the matrices (\ref{L_tri&L*_diag}). 
Therefore the elements (\ref{LP:Ln(0)}) give a Leonard pair on $L_n^{(0)}$ which is isomorphic to 
the Leonard pair $L, L^*$. 
\end{exam}

\begin{exam}
\label{exam:LP_Ln(1)}
Assume that $n\geq 1$ is an odd integer. Set the parameters 
$$
(r,s,d)=\left(\frac{1}{2},-\frac{1}{2},\frac{n-1}{2}\right).
$$
The settings of $r$ and $s$ fit the hypothesis of Theorem \ref{thm:LP}. 
Using Lemma \ref{lem:Ln(1)} yields that the matrices representing
\begin{gather}
\label{LP:Ln(1)}
\frac{E^2+F^2+\Lambda-1}{4}-\frac{H^2}{8}, \qquad \frac{n-H}{4}-\frac{1}{2}
\end{gather}
with respect to the ordered basis $\{v_i^{(1)}\}_{i=0}^d$ for $L_n^{(1)}$ are exactly the matrices (\ref{L_tri&L*_diag}). 
Therefore the elements (\ref{LP:Ln(1)}) give a Leonard pair on $L_n^{(1)}$ which is isomorphic to 
the Leonard pair $L, L^*$. 
%
\end{exam}

Fix a vertex $x$ of $\frac{1}{2}H(D,2)$. Let $\T(x)$ denote the Terwilliger algebra of $\frac{1}{2}H(D,2)$ with respect to $x$. 
By \cite[Theorem 6.4]{halved:2023} each $\T(x)$-module is a $\U_e$-module. By \cite[Corollary 6.6]{halved:2023} the irreducible $\U_e$-modules 
\begin{align*}
&L_{D-2k}^{(0)} 
\qquad \hbox{for all even integers $k$ with $0\leq k\leq \lfloor\frac{D}{2}\rfloor$};
\\
&L_{D-2k}^{(1)} 
\qquad \hbox{for all odd integers $k$ with $0\leq k\leq \lfloor\frac{D-1}{2}\rfloor$}
\end{align*}
are all irreducible $\T(x)$-modules up to isomorphism.  
Let $\mathbf A$ be the adjacency matrix of $\frac{1}{2}H(D,2)$. Let $\mathbf A^*(x)$ be the dual adjacency matrix of $\frac{1}{2}H(D,2)$ with respect to $x$.
The generators $\mathbf A$ and $\mathbf A^*(x)$ of $\T(x)$ act on the above irreducible $\U_e$-modules as 
$$
\frac{E^2+F^2+\Lambda-D}{2}-\frac{H^2}{4},
\qquad 
H,
$$
respectively. 
In the case of odd $D$, Examples \ref{exam:LP_Ln(0)} and \ref{exam:LP_Ln(1)} illustrate how the irreducible $\T(x)$-modules are related to Theorem \ref{thm:LP}; combined with Corollary \ref{cor:almost}, the Leonard pair $\mathbf A, \mathbf A^*(x)$ on any irreducible $\T(x)$-module is dual almost bipartite.

\subsection*{Funding statement}
The research was supported by the National Science and Technology Council of Taiwan under the project NSTC 114-2115-M-008-011.

\subsection*{Acknowledgements}
The author thanks Paul Terwilliger and Kazumasa Nomura for their valuable comments on the manuscript.

\subsection*{Conflict-of-interest statement}
The author has no conflict of interest to declare.

\subsection*{Data availability}

No data was used for the research described in the article.

\bibliographystyle{amsplain}
\bibliography{MP}

\providecommand{\bysame}{\leavevmode\hbox to3em{\hrulefill}\thinspace}
\providecommand{\MR}{\relax\ifhmode\unskip\space\fi MR }
\providecommand{\MRhref}[2]{%
  \href{http://www.ams.org/mathscinet-getitem?mr=#1}{#2}
}
\providecommand{\href}[2]{#2}
\begin{thebibliography}{10}

\bibitem{BannaiIto1984}
E.~Bannai and T.~Ito, \emph{{Algebraic Combinatorics I: Association Schemes}},
  Benjamin-Cummings, Menlo Park, 1984.

\bibitem{dualpolar:2022}
P.-A. Bernard, N.~Cramp\'{e}, and L.~Vinet, \emph{{The Terwilliger algebra of
  symplectic dual polar graphs, the subspace lattices and
  $U_q(\mathfrak{sl}_2)$}}, Discrete Mathematics \textbf{345} (2022), 113169.

\bibitem{Johnson:2021}
\bysame, \emph{{Entanglement of free fermions on Johnson graphs}}, Journal of
  Mathematical Physics \textbf{64} (2023), 061903.

\bibitem{SH:2020}
S.~Bockting-Conrad and H.-W. Huang, \emph{{The universal enveloping algebra of
  $\mathfrak{sl}_2$ and the Racah algebra}}, Communications in Algebra
  \textbf{48} (2020), 1022--1040.

\bibitem{DRG_book1989}
A.E. Brouwer, A.M. Cohen, and A.~Neumaier, \emph{{Distance-Regular Graphs}},
  Springer-Verlag, Berlin, 1989.

\bibitem{Brown2013}
G.M.F. Brown, \emph{{Totally bipartite/abipartite Leonard pairs and Leonard
  triples of Bannai/Ito type}}, Electronic Journal of Linear Algebra
  \textbf{26} (2013), 258--299.

\bibitem{Hahn:2019-2}
N.~Cramp\'{e}, E.~Ragoucy, L.~Vinet, and A.~Zhedanov, \emph{{Truncation of the
  reflection algebra and the Hahn algebra}}, Journal of Physics A: Mathematical
  and Theoretical \textbf{52} (2019), 35LT01.

\bibitem{Hahn:2019-1}
L.~Frappat, J.~Gaboriaud, L.~Vinet, S.~Vinet, and A.~Zhedanov, \emph{{The Higgs
  and Hahn algebras from a Howe duality perspective}}, Physics Letters A
  \textbf{383} (2019), 1531--1535.

\bibitem{gvz2014}
V.X. Genest, L.~Vinet, and A.~Zhedanov, \emph{{Superintegrability in two
  dimensions and the Racah--Wilson algebra}}, Letters in Mathematical Physics
  \textbf{104} (2014), 931--952.

\bibitem{gvz2013}
\bysame, \emph{{The equitable Racah algebra from three $\mathfrak{su}(1, 1)$
  algebras}}, Journal of Physics A: Mathematical and Theoretical \textbf{47}
  (2014), 025203 (12 pages).

\bibitem{R&LD2014}
\bysame, \emph{{The Racah algebra and superintegrable models}}, Journal of
  Physics: Conference Series \textbf{512} (2014), 012011 (15 pages).

\bibitem{hypercube2002}
J.T. Go, \emph{{The Terwilliger algebra of the hypercube}}, European Journal of
  Combinatorics \textbf{23} (2002), 399--429.

\bibitem{Huang:CG&Hamming}
H.-W. Huang, \emph{{The Clebsch--Gordan rule for $U(\mathfrak{sl}_2)$, the
  Krawtchouk algebras and the Hamming graphs}}, SIGMA \textbf{19} (2023), 017,
  19 pages.

\bibitem{Huang:CG&Johnson}
\bysame, \emph{{The Clebsch--Gordan coefficients of $U(\mathfrak{sl}_2)$ and
  the Terwilliger algebras of Johnson graphs}}, Journal of Combinatorial
  Theory, Series A \textbf{203} (2024), 105833.

\bibitem{Huang:CG&Grassmann}
\bysame, \emph{{An imperceptible connection between the Clebsch--Gordan
  coefficients of $U_q(\mathfrak{sl}_2)$ and the Terwilliger algebras of
  Grassmann graphs}}, Journal of Combinatorial Theory, Series A \textbf{214}
  (2025), 106028.

\bibitem{SH:2019-1}
H.-W. Huang and S.~Bockting-Conrad, \emph{{Finite-dimensional irreducible
  modules of the Racah algebra at characteristic zero}}, SIGMA \textbf{16}
  (2020), 018, 17 pages.

\bibitem{SH:2017-1}
\bysame, \emph{{The Casimir elements of the Racah algebra}}, Journal of Algebra
  and Its Applications \textbf{922} (2021), 2150135 (22 pages).

\bibitem{halved:2023}
H.-W. Huang and C.-Y. Wen, \emph{{A connection behind the Terwilliger algebras
  of $H(D,2)$ and $\frac{1}{2}H(D,2)$}}, Journal of Algebra \textbf{634}
  (2023), 456--479.

\bibitem{halved:2024}
\bysame, \emph{{Johnson graphs as slices of a hypercube and an algebra
  homomorphism from the universal Racah algebra into $U(\mathfrak{sl}_2)$}},
  Journal of Algebraic Combinatroics \textbf{61} (2025), article number 56.

\bibitem{Hamming:2015}
M.A. Jafarizadeh, S.~Nami, and F.~Eghbalifam, \emph{{Entanglement entropy in
  the Hamming networks}}, {arXiv:1503.04986}.

\bibitem{Koe2010}
R.~Koekoek, P.~Lesky, and R.~Swarttouw, \emph{Hypergeometric orthogonal
  polynomials and their $q$-analogues}, Springer Monographs in Mathematics,
  Springer-Verlag, Berlin, 2010.

\bibitem{Johnson:2007}
F.~Levstein and C.~Maldonado, \emph{{The Terwilliger algebra of the Johnson
  schemes}}, Discrete Mathematics \textbf{307} (2007), 1621--1635.

\bibitem{Masuda2024}
S.~Masuda, \emph{{Doubly almost bipartite Leonard pairs}}, Ph.D. thesis,
  {Portland State University}, 2024.

\bibitem{Doob2021}
J.V.S. Morales, \emph{{On quantum adjacency algebras of Doob graphs and their
  irreducible modules}}, Journal of Algebraic Combinatorics \textbf{54} (2021),
  979--998.

\bibitem{Doob2023}
\bysame, \emph{{On standard bases of irreducible modules of Terwilliger
  algebras of Doob schemes}}, Journal of Algebraic Combinatorics \textbf{58}
  (2023), 913--931.

\bibitem{nomura2025-2}
K.~Nomura and P.~Terwilliger, \emph{{The Norton-balanced condition for
  $Q$-polynomial distance-regular graphs}}, The Electronic Journal of
  Combinatorics \textbf{32} (2025), \#P1.49.

\bibitem{bipart&so6}
P.~Terwilliger, \emph{{$2$-Homogeneous bipartite distance-regular graphs and
  the quantum group $U_q'(\mathfrak{so}_6)$}}, arXiv:2506.02190.

\bibitem{TerAlgebraI}
\bysame, \emph{{The subconstituent algebra of an association scheme (part I)}},
  Journal of Algebraic Combinatorics \textbf{1} (1992), 363--388.

\bibitem{TerAlgebraII}
\bysame, \emph{{The subconstituent algebra of an association scheme (part
  II)}}, Journal of Algebraic Combinatorics \textbf{2} (1993), 73--103.

\bibitem{TerAlgebraIII}
\bysame, \emph{{The subconstituent algebra of an association scheme (part
  III)}}, Journal of Algebraic Combinatorics \textbf{2} (1993), 177--210.

\bibitem{Askeyscheme}
\bysame, \emph{{An algebraic approach to the Askey scheme of orthogonal
  polynomials}}, {Orthogonal polynomials and special functions: Computation and
  applications} (Berlin) (F.~Marcell\'{a}n and W.~Van Assche, eds.), Lecture
  Notes in Mathematics 1883, Springer, 2006, pp.~255--330.

\bibitem{Wen2024}
C.-Y. Wen, \emph{{A study of relationships between the universal Hahn algebra
  and the universal Racah algebra with $U(\mathfrak{sl}_2)$}}, Ph.D. thesis,
  {National Yang Ming Chiao Tung University}, 2024.

\end{thebibliography}

\end{document}